\documentclass[11pt]{article}
\usepackage{amsmath,amssymb,amsbsy,amsfonts,amsthm,
   latexsym,amsopn,amstext,amsxtra,euscript,amscd}

\begin{document}

\newtheorem{prop}{Proposition}[section]
\newtheorem{lem}[prop]{Lemma}
\newtheorem{cor}[prop]{Corollary}
\newtheorem{conj}[prop]{Conjecture}
\newtheorem{defi}[prop]{Definition}
\newtheorem{thm}[prop]{Theorem}
\newtheorem{rem}[prop]{Remark}
\newtheorem{rems}[prop]{Remarks}
\newtheorem{fac}[prop]{Fact}
\newtheorem{facs}[prop]{Facts}
\newtheorem{com}[prop]{Comments}
\newtheorem{prob}{Problem}
\newtheorem{problem}[prob]{Problem}
\newtheorem{ques}{Question}
\newtheorem{question}[ques]{Question}


\def\scr{\scriptstyle}
\def\\{\cr}
\def\({\left(}
\def\){\right)}
\def\[{\left[}
\def\]{\right]}
\def\<{\langle}
\def\>{\rangle}
\def\fl#1{\left\lfloor#1\right\rfloor}
\def\rf#1{\left\lceil#1\right\rceil}
\def\defn{\noindent{\bf Definition\/}\ \ }
\def\rem{\noindent{\bf Remarks\/}\ \ }
\def\card#1{\vphantom{#1}^{\#}}

\def\Z{\mathbb Z}
\def\R{\mathbb R}
\def\N{\mathbb N}
\def\cS{\mathcal S}
\def\e{{\bf e}}
\def\cR{\mathcal R}
\def\cP{\mathcal P}
\def\cQ{\mathcal Q}
\def\cT{\mathcal T}
\def\rad{{\mathrm{rad}\/}}
\def\eps{\varepsilon}
\def\ep{{\mathbf{\,e}}_p}
\def\epp{{\mathbf{\,e}}_{p-1}}

\def\xxx{\vskip5pt\hrule\vskip5pt}
\def\yyy{\vskip5pt\hrule\vskip2pt\hrule\vskip5pt}

\newcommand{\comm}[1]{\marginpar {\fbox{#1}}}
\newcommand{\1}{1\!{\mathrm l}}
\newcommand{\croix}{
set\frown}




\newfont{\teneufm}{eufm10}
\newfont{\seveneufm}{eufm7}
\newfont{\fiveeufm}{eufm5}
%
%
\newfam\eufmfam
                  \textfont\eufmfam=\teneufm \scriptfont\eufmfam=\seveneufm
                  \scriptscriptfont\eufmfam=\fiveeufm

%
%
\def\frak#1{{\fam\eufmfam\relax#1}}
%


\def\bbbr{{\rm I\!R}} 
\def\bbbc{{\rm I\!C}} 
\def\bbbm{{\rm I\!M}}
\def\bbbn{{\rm I\!N}} 
\def\bbbf{{\rm I\!F}}
\def\bbbh{{\rm I\!H}}
\def\bbbk{{\rm I\!K}}
\def\bbbl{{\rm I\!L}}
\def\bbbp{{\rm I\!P}}
\newcommand{\lcm}{{\rm lcm}}
\def\bbbone{{\mathchoice {\rm 1\mskip-4mu l} {\rm 1\mskip-4mu l}
{\rm 1\mskip-4.5mu l} {\rm 1\mskip-5mu l}}}
\def\bbbc{{\mathchoice {\setbox0=\hbox{$\displaystyle\rm C$}\hbox{\hbox
to0pt{\kern0.4\wd0\vrule height0.9\ht0\hss}\box0}}
{\setbox0=\hbox{$\textstyle\rm C$}\hbox{\hbox
to0pt{\kern0.4\wd0\vrule height0.9\ht0\hss}\box0}}
{\setbox0=\hbox{$\scriptstyle\rm C$}\hbox{\hbox
to0pt{\kern0.4\wd0\vrule height0.9\ht0\hss}\box0}}
{\setbox0=\hbox{$\scriptscriptstyle\rm C$}\hbox{\hbox
to0pt{\kern0.4\wd0\vrule height0.9\ht0\hss}\box0}}}}
\def\bbbq{{\mathchoice {\setbox0=\hbox{$\displaystyle\rm
Q$}\hbox{\raise
0.15\ht0\hbox to0pt{\kern0.4\wd0\vrule height0.8\ht0\hss}\box0}}
{\setbox0=\hbox{$\textstyle\rm Q$}\hbox{\raise
0.15\ht0\hbox to0pt{\kern0.4\wd0\vrule height0.8\ht0\hss}\box0}}
{\setbox0=\hbox{$\scriptstyle\rm Q$}\hbox{\raise
0.15\ht0\hbox to0pt{\kern0.4\wd0\vrule height0.7\ht0\hss}\box0}}
{\setbox0=\hbox{$\scriptscriptstyle\rm Q$}\hbox{\raise
0.15\ht0\hbox to0pt{\kern0.4\wd0\vrule height0.7\ht0\hss}\box0}}}}
\def\bbbt{{\mathchoice {\setbox0=\hbox{$\displaystyle\rm
T$}\hbox{\hbox to0pt{\kern0.3\wd0\vrule height0.9\ht0\hss}\box0}}
{\setbox0=\hbox{$\textstyle\rm T$}\hbox{\hbox
to0pt{\kern0.3\wd0\vrule height0.9\ht0\hss}\box0}}
{\setbox0=\hbox{$\scriptstyle\rm T$}\hbox{\hbox
to0pt{\kern0.3\wd0\vrule height0.9\ht0\hss}\box0}}
{\setbox0=\hbox{$\scriptscriptstyle\rm T$}\hbox{\hbox
to0pt{\kern0.3\wd0\vrule height0.9\ht0\hss}\box0}}}}
\def\bbbs{{\mathchoice
{\setbox0=\hbox{$\displaystyle     \rm S$}\hbox{\raise0.5\ht0\hbox
to0pt{\kern0.35\wd0\vrule height0.45\ht0\hss}\hbox
to0pt{\kern0.55\wd0\vrule height0.5\ht0\hss}\box0}}
{\setbox0=\hbox{$\textstyle        \rm S$}\hbox{\raise0.5\ht0\hbox
to0pt{\kern0.35\wd0\vrule height0.45\ht0\hss}\hbox
to0pt{\kern0.55\wd0\vrule height0.5\ht0\hss}\box0}}
{\setbox0=\hbox{$\scriptstyle      \rm S$}\hbox{\raise0.5\ht0\hbox
to0pt{\kern0.35\wd0\vrule height0.45\ht0\hss}\raise0.05\ht0\hbox
to0pt{\kern0.5\wd0\vrule height0.45\ht0\hss}\box0}}
{\setbox0=\hbox{$\scriptscriptstyle\rm S$}\hbox{\raise0.5\ht0\hbox
to0pt{\kern0.4\wd0\vrule height0.45\ht0\hss}\raise0.05\ht0\hbox
to0pt{\kern0.55\wd0\vrule height0.45\ht0\hss}\box0}}}}
\def\bbbz{{\mathchoice {\hbox{$\sf\textstyle Z\kern-0.4em Z$}}
{\hbox{$\sf\textstyle Z\kern-0.4em Z$}}
{\hbox{$\sf\scriptstyle Z\kern-0.3em Z$}}
{\hbox{$\sf\scriptscriptstyle Z\kern-0.2em Z$}}}}
\def\ts{\thinspace}

\def\squareforqed{\hbox{\rlap{$\sqcap$}$\sqcup$}}
\def\qed{\ifmmode\squareforqed\else{\unskip\nobreak\hfil
\penalty50\hskip1em\null\nobreak\hfil\squareforqed
\parfillskip=0pt\finalhyphendemerits=0\endgraf}\fi}

\def\cA{{\mathcal A}}
\def\cB{{\mathcal B}}
\def\cC{{\mathcal C}}
\def\cD{{\mathcal D}}
\def\cE{{\mathcal E}}
\def\cF{{\mathcal F}}
\def\cG{{\mathcal G}}
\def\cH{{\mathcal H}}
\def\cI{{\mathcal I}}
\def\cJ{{\mathcal J}}
\def\cK{{\mathcal K}}
\def\cL{{\mathcal L}}
\def\cM{{\mathcal M}}
\def\cN{{\mathcal N}}
\def\cO{{\mathcal O}}
\def\cP{{\mathcal P}}
\def\cQ{{\mathcal Q}}
\def\cR{{\mathcal R}}
\def\cS{{\mathcal S}}
\def\cT{{\mathcal T}}
\def\cU{{\mathcal U}}
\def\cV{{\mathcal V}}
\def\cW{{\mathcal W}}
\def\cX{{\mathcal X}}
\def\cY{{\mathcal Y}}
\def\cZ{{\mathcal Z}}

\def\ord{{\mathrm{ord}}}
\def\Nm{{\mathrm{Nm}}}
\renewcommand{\vec}[1]{\mathbf{#1}}

\def \F{{\bbbf}}
\def \L{{\bbbl}}
\def \K{{\bbbk}}
\def \C{{\bbbc}}
\def \Z{{\bbbz}}
\def \N{{\bbbn}}
\def \Q{{\bbbq}}
\def\E{{\mathbf E}}
\def\H{{\mathbf H}}
\def\G{{\mathcal G}}
\def\O{{\mathcal O}}
\def\cS{{\mathcal S}}
\def \R{{\bbbr}}
\def\Fp{\F_p}
\def \fp{\Fp^*}
\def\\{\cr}
\def\({\left(}
\def\){\right)}
\def\fl#1{\left\lfloor#1\right\rfloor}
\def\rf#1{\left\lceil#1\right\rceil}

\def\Zm{\Z_m}
\def\Zt{\Z_t}
\def\Zp{\Z_p}
\def\Um{\cU_m}
\def\Ut{\cU_t}
\def\Up{\cU_p}

\def\ep{{\mathbf{e}}_p}

\def \Prob{{\mathrm {}}}





\title{\bf On the period of the continued fraction expansion 
of ${\sqrt {2^{2n+1}+1}}$}

\author{
{\sc Yann Bugeaud} \\
{Universit\'e Louis Pasteur}\\
{UFR de math\'ematiques} \\
{7 rue Ren\'e Descartes, 67084 Strasbourg, France} \\
{\tt bugeaud@math.u-strasbg.fr}
\and
{\sc Florian~Luca} \\
{Instituto de Matem{\'a}ticas}\\
{ Universidad Nacional Aut\'onoma de M{\'e}xico} \\
{C.P. 58180, Morelia, Michoac{\'a}n, M{\'e}xico} \\
{\tt fluca@matmor.unam.mx}}

\date{\today}

\maketitle

\begin{abstract}

\end{abstract}
In this paper, we prove that the period of the continued fraction 
expansion of ${\sqrt {2^{n}+1}}$ tends to infinity when $n$ tends to infinity 
through odd positive integers.

\section{Introduction}

It is, in general, very hard to predict the features of the continued 
fraction expansion of a given positive real number. If the number in question 
is of the form ${\sqrt {d}}$, where $d$ is a positive integer which 
is not a square, then its continued fraction expansion is of the form 
$[a_0,\{a_1,\dots,a_{r-1}, 2a_0\}]$, where we use 
$\{\dots\}$ to emphasize the period of the expansion. It is known that 
$a_1,\dots,a_{r-1}$ is a palindrome; i.e., $a_i=a_{r-i}$ holds 
for all $i=1,\dots,r-1$. The lenght $r$ of the period is at least $1$ 
(and this is achieved, for example, for square free numbers 
$d$ of the form $k^2+1$ with some positive integer $k$), and  
$r\ll {\sqrt {d}}\log d$ (see ~\cite{Hua}). Here, and in all what 
follows, we use the Vinogradov symbols $\ll$ and $\gg$, 
as well as the Landau symbols $O$ and $o$, with their usual meanings.

It is believed that for ``most'' 
$d$ the above upper bound is close to the truth. 
For this and other open problems concerning the behavior of $r$ as 
a function of $d$, we refer the reader to Lenstra's paper ~\cite{Len}.

When $d$ is restricted to run through certain parametrized families, 
occasionally some regular patterns appear. For example, Schinzel 
(see ~\cite{Sch1}, ~\cite{Sch2}) proved that if $f(X)$ is a 
non constant polynomial with integer 
coefficients and positive leading term 
satisfying certain assumptions (for example, of odd degree, 
or of even degree but whose leading term is not a square of a positive 
integer), then the length of the continued fraction expansions 
of ${\sqrt {f(n)}}$ can become arbitrarily large as $n$ goes to infinity. 

In this paper, we look at a problem similar to Schinzel's problem 
mentioned above when the polynomial $f(n)$ is replaced by a power sum 
over ${\Z}$ satisfying suitable assumptions. That is, 
let $\ell\ge 1$, $a_i$ and $b_i$ be non zero integers for $i=1,\dots,\ell$, 
with $b_1>b_2>\dots>b_\ell\ge 1$, and set
\begin{equation}
\label{eq:f}
f(n)=\sum_{i=1}^\ell a_i b_i^n.
\end{equation}
We call $b_1,\dots,b_\ell$ the {\it roots} 
of the form $f(n)$ and $a_1,\dots,a_\ell$ 
its {\it coefficients}.
To follow
standard notations (see, ~\cite{CZ1}, for example), we write 
${\cal E}_{\Z}$ for the ring of all such forms together with 
the constant $0$ form. If $R$ is any subring of $\C$, we write 
$R{\cal E}_{\Z}$ for the ring $R\otimes_{\Z} {\cal E}_{\Z}$, 
which is the ring of power sums $f(n)$ given by 
formula ~\eqref{eq:f}, but where the coefficients $a_i$ are allowed 
to be in $R$. As usual, we write ${\overline {\Q}}$ for the field of 
algebraic numbers.
Whenever we write ${\sqrt {f(n)}}$ for some $f(n)\in {\Q}{\cal E}_{\Z}$, 
we implicitly mean that $a_1>0$. In this way, we ensure that 
the above square root is real for all but finitely many values of the positive 
integer $n$.
\medskip

{\bf Acknowledgments.}  This paper was written during a visit 
of Y.~B. at the Mathematical Institute of the UNAM in Morelia in 
January 2004. He thanks this Institute for its hospitality. Both 
authors thank Pietro Corvaja and Umberto Zannier for a copy of 
~\cite{CZ2}. Both 
authors were supported in part by the joint Project France-Mexico 
ANUIES-ECOS M01-M02.

\section{Results}

In order to prove our main result, we shall assume that 
our form $f(n)\in {\Q}{\cal E}_{\Z}$ satisfies the following condition:

\medskip

{\bf Hypothesis (H).} 
{\it There do not exist an integer $j\in \{0,1\}$, a number 
$\delta<1/2$, and forms $g(n)$ and $h(n)$ in ${\Q}{\cal E}_{\Z}$, 
such that both the relation
$$f(2n+j)=h(n)^2+g(n)$$and the estimate
$$|g(n)|\ll |f(n)|^{\delta}$$hold for all positive integers $n$.}

\medskip

In this paper, we prove the following result.

\begin{thm}
\label{thm:Main}
Assume that $f(n)\in {\cal E}_{\Z}$ satisfies Hypothesis (H). 
Then ${\sqrt {f(n)}}$ is a rational number for at 
most finitely many positive 
integers $n$. Moreover, the length $r(n)$ of the period of the continued 
fraction expansion of ${\sqrt {f(n)}}$ tends to infinity with $n$.
\end{thm}

It is likely that Theorem ~\ref{thm:Main} 
remains true even for certain forms $f(n)\in 
{\cal E}_{\Z}$ (or ${\Q}{\cal E}_{\Z}$) 
which do not satisfy the above Hypothesis (H). However, 
note that some restrictions must be imposed as, for example, 
${\sqrt {h(n)^2+1}}=[h(n),\{2h(n)\}]$ holds for all forms 
$h(n)\in {\cal E}_{\Z}$ whose coefficients $a_i$ are positive 
for $i=1,\dots,\ell$, while the example $f(n)=h(n)^2$ with 
$h(n)\in {\Q}{\cal E}_{\Z}$ shows that 
${\sqrt {f(n)}}$ can be a rational number with a bounded 
denominator for all positive integers $n$. See Section 5 for further 
remarks.

While the above Hypothesis (H) seems cumbersome to verify, we note that 
it trivially holds if none of the two positive integers $a_1$ or $a_1b_1$
is a square. In particular, Theorem ~\ref{thm:Main} applies to the form 
$f(n)=2\cdot 4^n+1$  mentioned in the title of the present paper.

We also point out that Theorem ~\ref{thm:Main} 
gives a partial answer to a problem 
specifically raised at the end of ~\cite{CZ2}, where it is predicted 
that the period of the continued fraction expansion of 
${\sqrt {f(n)}}$ tends to infinity with $n$ once 
$f(n)\in {\Q}{\cal E}_{\Z}$ satisfies certain ``suitable assumptions'', 
which is the case here. 

As predicted in the concluding remarks of 
~\cite{CZ2}, the proof of Theorem ~\ref{thm:Main} 
uses the Subspace Theorem, much 
in the spirit of the papers ~\cite{CZ1} and ~\cite{CZ2}.

\section{Preparations}

In this section, we review some standard notions of algebraic number theory
(see, for example, \cite{Cohn,Nark,StTal}), and
Diophantine approximations.

Let $\L$ be an algebraic number field of degree $D$ over $\Q$.
Denote its ring of integers by $O_{\L}$ and its collection of places
by $\cM_{\L}$. For a fractional ideal $\cI$ of $\L$, we denote by
$\Nm_{\L}(\cI)$ its norm. We recall that $\Nm_{\L}(\cI)=\#(O_{\L}/{\cI})$
if $\cI$ is an ideal of $O_{\L}$, and the norm map
is extended multiplicatively
(using unique factorization) to all the fractional ideals
of $\L$.

For a prime  ideal $\cP$, we denote by $\ord_{\cP} (x)$
the order at which it appears in the factorization of the
principal ideal $[x]$ generated by $x$ inside $\L$.

For $\mu\in \cM_{\L}$ and $x\in \L$, we define the absolute value
$|x|_{\mu}$ as follows:
\begin{description}
\item[(i)] $|x|_{\mu}=|\sigma(x)|^{1/D}$
if $\mu$ corresponds to the
embedding $\sigma: \L\mapsto \R$;

\item[(ii)] $|x|_{\mu}=|\sigma(x)|^{2/D}=|{\overline{\sigma}}(x)|^{2/D}$
if $\mu$ corresponds to the pair of complex
conjugate embeddings $\sigma,{\overline{\sigma}}:\L\mapsto
\bbbc$;

\item[(iii)] $|x|_{\mu}=\Nm_{\L}(\cP)^{\ord_{\cP}(x)}$ if $\mu$ corresponds
to the nonzero prime ideal ${\cP}$ of $O_{\L}$.
\end{description}

In  case~(i) or~(ii) we say that $\mu$ is {\it real infinite} or
{\it complex infinite},
respectively; in case~(iii) we say that that $\mu$ is {\it finite}.

These absolute values satisfy the {\it product formula\/}
$$
\prod_{\mu\in {\cM}_{\L}} |x|_{\mu}=1,\qquad  \text{ for all } x\in \L^*.
$$

Our basic tool is the following simplified version 
of a result  of Schlickewei (see ~\cite{Schm1}, ~\cite{Schm2}), 
which is commonly known as the Subspace Theorem.

\begin{lem}
\label{lem:S-unit}
Let $\L$ be an algebraic number field of degree $D$. 
Let $\cS$ be a finite set of places of $\L$ containing all
the infinite ones. Let  $\{L_{1,\mu},\ldots,L_{M,\mu}\}$ for $\mu\in \cS$ be
linearly  independent sets of linear forms in $M$ variables 
with coefficients in $\L$.
Then, for every fixed $0<\varepsilon<1$,
the set $\cX$
of solutions $\vec{x}=(x_1,\dots,x_M) \in \Z^M\backslash \{0\}$ to the
inequality
\begin{equation}
\label{eq:fundineq}
\prod_{\mu\in \cS}\prod_{i=1}^M
|L_{i,\mu}(\vec{x})|_{\mu} <{\rm max}\{|x_i|~|~i=1,\dots,M\}^{-\varepsilon},
\end{equation}
is contained in finitely many proper linear subspaces of $\Q^M$. 
\end{lem}

\section{Proofs}

Throughout this section, $C_1,C_2,\dots$ are effectively computable constants
which are either absolute, or 
depend on the given data (usually, a form $f(n)\in {\cal E}_{\Z}$). 

The following result is a variation of Lemma 1 from ~\cite{CZ1}.

\begin{lem}
\label{lem:CZ1}
There exists an absolute constant 
$C_1$ such that the following holds.  
If $b$ is any positive integer and $f(n)\in {\cal E}_{\Z}$ 
(not necessarily satisfying Hypothesis (H)) are such that 
for infinitely many positive integers 
$n$ the denominator of the rational number $f(n)/b^n$ is less than 
$\exp(C_1n)$, then $b~|~b_i$ for all $i=1,\dots, \ell$.
\end{lem}

\begin{proof} We shall choose $C_1=\log 2/2$. Without any loss of 
generality, we may assume that ${\gcd}(b_1,\dots,b_\ell)=1$. We then 
have to prove that $b=1$. Assume that this is not so, and assume 
further that $b$ is prime (if not, we replace $b$ by a prime factor of 
it). Finally, it is clear that 
we may assume that none of the roots of $f(n)$ is a multiple of $b$, for if 
not, we may replace $f(n)$ by 
$$
\sum_{\substack{{1\le i\le \ell}\\ {b~{\not|}~b_i}}} a_ib_i^n.
$$
We now apply Lemma ~\ref{lem:S-unit} as in the proof 
of Lemma 1 in ~\cite{CZ1}. We let $\L=\Q$, $M=\ell$, and  
${\cal S}$ be the set of places of $\L$ consisting of $\infty$, $b$, and all 
prime factors of $b_i$ for $i=1,\dots, \ell$. For $\mu\in {\cal S}\backslash 
\{b\}$ and a vector ${\bf x}=(x_1,\dots,x_M)$ we put 
$L_{i,\mu}({\bf x})=x_i$
for $i=1,\dots, M$, while for $\mu=b$ we put 
$L_{1,b}({\bf x})=\sum_{i=1}^{M} a_ix_i$ and $L_{i,b}({\bf x})=x_i$ 
for $i=2,\dots,M$. We evaluate the double product appearing in 
the statement of Lemma ~\ref{lem:S-unit} for 
${\bf x}=(b_1^n,\dots,b_{\ell}^n)$. We note that $x_i$ are integers 
for all $i=1,\dots,M$. The calculation from page 322 in ~\cite{CZ1} 
shows that
\begin{equation}
\label{ineq:DP}
\prod_{\mu\in {\cal S}} \prod_{i=1}^{M} |L_{i,\mu}({\bf x})|_{\mu}=
|L_{1,b}({\bf x})|_b\le b^{-n}\cdot 2^{n/2}\le b^{-n/2}=(b_1^n)^
{-\varepsilon_0},
\end{equation}
where $\varepsilon_0=\log p/(2\log b_1)$. Since 
$b_1^n={\max}\{|x_i|~|~i=1,\dots,M\}$, it follows easily that 
the above inequality ~\eqref{ineq:DP} implies that our points 
${\bf x}$ and linear forms $L_{i,\mu}$ for $i=1,\dots, M$, and 
$\mu\in {\cal S}$ fulfill inequality ~\eqref{eq:fundineq} with 
$\varepsilon=\varepsilon_0$. Now Lemma  ~\ref{lem:S-unit} asserts 
that there are finitely many proper subspaces of ${\Q}^{M}$ 
of equations of the form $\sum_{i=1}^{M} c_ix_i=0$ 
with $c_i\in \Q$ for $i=1,\dots,M$, not all zero, 
such that every point ${\bf x}\in \Z^{M}$ satisfying the above inequality 
~\eqref{ineq:DP} lies on one of these subspaces. This in turns 
gives us equations of the form 
\begin{equation}
\label{eq:LRS=zero}
\sum_{i=1}^{M} c_ib_i^n=0.
\end{equation}
Since each one of the above equations gives the set of zeros of
a linear recurrent sequence having a dominant root (note that 
at least one of the coefficients $c_i$ is non zero), it follows that
each one of these equations can  
have only finitely many positive integer solutions $n$. 
\end{proof}

Let $f(n)\in {\Q}{\cal E}_{\Z}$ be some form, not necessarily satisfying
Hypothesis (H). Replacing $f(n)$ by $f(2n+j)$ for $j=\{0,1\}$, it follows 
that we may replace $b_i$ by $b_i^2$ and $a_i$ by $a_ib_i^j$ for 
$i=1,\dots,\ell$. In particular, we may assume that $b_1$ is a square.

\begin{lem}
\label{lem:fundunit}
Let $f(n)\in {\Q}{\cal E}_{\Z}$. Assume that $f(n)$ satisfies 
Hypothesis (H). Then there exists a computable positive constant 
$C_2$, depending only on $f(n)$, such that if $C$ is any fixed constant 
and if $(X(n),Y(n))$ are positive  integers such that the inequality
$$
|X(n)^2-f(n)Y(n)^2|<C
$$
holds, then $Y(n)>\exp(C_2n)$ holds for all positive integers $n$ with 
finitely many exceptions. 
\end{lem}

\begin{proof}
We write $f(n)=a_1b_1^n(1+\delta(n)),$ where
$$
\delta(n)=\sum_{i=2}^{\ell} \frac{a_i}{a_1}\(\frac{b_i}{b_1}\)^n.
$$
Note that $\delta(n)=0$ in ${\Q}{\cal E}_{\Z}$ 
if and only if $\ell=1$. If $\ell\ge 2$, 
we then let $\beta=b_1/b_2$, and observe that $\beta>1$ 
and that $\delta(n)=O(\beta^{-n})$. 
We let $k$ be a positive integer such that $\beta^k>b_1$. 
Clearly, we can choose $k=\lfloor \log b_1/\log\beta\rfloor +1$. 
Writing $\alpha={\sqrt {a_1}}$, we note that we have the approximation
\begin{eqnarray*}
{\sqrt {f(n)}} & = & \alpha b_1^{n/2}{\sqrt {1+\delta(n)}}\\
& = & \alpha b_1^{n/2}\(\sum_{i=0}^k \binom{1/2}{i} \delta(n)^{i} +
O(\delta(n)^{k+1})\)\\
& = & \alpha b_1^{n/2}\sum_{i=0}^k \binom{1/2}{i} \delta(n)^{i}+
O\(b_1^{-n/2}\beta^{-n}\).
\end{eqnarray*}
Note that 
$$
\alpha b_1^{n/2} \sum_{i=0}^k \binom{1/2}{i} \delta(n)^{i}=\alpha \cdot 
\frac{f_1(n)}{b_1^{(k-1/2)n}},
$$
where $f_1(n)\in {\Q}{\cal E}_{\Z}.$ Thus, we may write that 
$$
{\sqrt {f(n)}}=\alpha \frac{f_1(n)}{b_1^{(k-1/2)n}}+O(b_1^{-n/2}\beta^{-n}),
$$
where we take $f_1(n)=0$ if $\ell=1$. Note also that all the prime factors 
of the roots of $f_1(n)$ are among the prime factors of the roots of 
$f(n)$. 
Assume now that $C$ is some fixed 
positive constant and that $(X(n),Y(n))$ is a pair 
of positive integers such that 
\begin{equation}
\label{eq:C}
|X(n)^2-f(n)Y(n)^2|<C.
\end{equation}
Then, since
$$
f(n)=\(\alpha \frac{f_1(n)}{b_1^{(k-1/2)n}}+O(b_1^{-n/2}
\beta^{-n})\)^2=\alpha^2 \(\frac{f_1(n)}{b_1^{(k-1/2)n}}\)^2+O\(\beta^{-n}\),
$$
we get that
$$
X(n)^2-f(n)Y(n)^2=X(n)^2-\alpha^2 \(\frac{f_1(n)}{b_1^{(k-1/2)n}}\)^2Y(n)^2+
O\(Y(n)^2\beta^{-n}\).
$$
We choose $C_2<\log \beta/2$, and infer that if $Y(n)<\exp(C_2n)$, then 
inequality ~\eqref{eq:C} leads to the conclusion that the inequality
$$
\left|X(n)^2-\alpha^2 \(\frac{f_1(n)}{b_1^{(k-1/2)n}}\)^2Y(n)^2\right|<2C
$$
holds for all but finitely many positive integers $n$. In turn, 
the above inequality implies that 
\begin{equation}
\label{eq:*}
\left|X(n)-\frac{\alpha f_1(n)}{b_1^{(k-1/2)n}}Y(n)
\right|\ll \frac{1}{b_1^{n/2} Y(n)}.
\end{equation}
The constant understood in $\ll$ above depends on $C$ and on the form 
$f(n)$. The above inequality ~\eqref{eq:*} is equivalent to 
\begin{equation}
\label{ineq:DP2}
\left|b_1^{(k-1/2)n}X(n)-\alpha f_1(n)Y(n)\right|
\ll \frac{b_1^{(k-1)n}} {Y(n)}.
\end{equation}
We now write 
$$
f_1(n)=\sum_{i=1}^{\ell'} a_i' (b_i')^n,
$$
where $b_1'>b_2'> \dots >b_{\ell'}\ge 1$. 
Note that $1\le \ell'\le 1+(\ell-1)+\dots+(\ell-1)^k$, and that 
$b_1'=b_1^{k}$ (recall that $b_1$ is a square). We are now all set to 
apply Lemma ~\ref{lem:S-unit}. We choose $\L=\Q[\alpha]$, $M=1+\ell'$, 
and ${\cal S}$ to be the set of all places of $\L$ 
(which is either $\Q$, or a real quadratic field, respectively) 
consisting of the infinite ones (one, or two of them, respectively), and 
the finite ones corresponding to primes in $\L$ lying above the prime 
factors of $b_1\dots b_{\ell}$. Note that all the prime 
factors of the $b_i'$'s are among the prime factors of the $b_i$'s. 
When $\mu\in {\cal S}$ is finite, 
we then put $L_{i,\mu}({\bf x})=x_i$ for $i=1,\dots,M$, 
while if $\mu$ is infinite
corresponding to the real embedding $\sigma:\L\mapsto \R$, we then 
put $L_{i,\mu}=x_i$ if $i\ne 2$, and 
$L_{i,\mu}=x_1-\sigma^{-1}(\alpha)(a_1'x_2+\dots+a'_{\ell'} x_{\ell'+1})$ 
if $i=2$. Note that if $x_i$ are rational integers (i.e., in $\Z$), then 
$|L_{i,\mu}({\bf x})|_{\mu}=
|x_1-\alpha(a_1'x_2+\dots+a'_{\ell'} x_{\ell'+1})|^{1/D}$
holds for all the infinite places $\mu\in {\cal S}$.
We now verify that if we take ${\bf x}=(x_1,\dots,x_M)$ 
as $x_1=b_1^{(k-1/2)n}X(n)$, and $x_i=(b_{i-1}')^n Y(n)$ for $i=2,\dots, M$, 
then inequality ~\eqref{ineq:DP2} implies that the inequality 
\begin{equation}
\label{eq:pell}
\prod_{\mu\in {\cal S}} \prod_{i=1}^M |L_{i,\mu}({\bf x})|_{\mu}\ll
 \frac{Y(n)^{M-1}}{b_1^{n/2}}
\end{equation}
holds. Observe that 
if $i\ne 2$, then 
$$
\prod_{\mu\in {\cal S}} |L_{i,\mu}({\bf x})|_{\mu}=\prod_{\mu\in {\cal S}}
|x_i|_{\mu},
$$
and by the product formula, the fact that $x_i\in \Z^*$ for all 
$i=1,\dots,M$, and the fact that ${\cal S}$ contains all infinite places 
and all the places corresponding to all the 
prime divisors of $b_i'$ for $i=1,\dots, \ell'$,
it follows easily that
\begin{equation}
\label{eq:prod1}
\prod_{\mu\in {\cal S}} |L_{1,\mu}({\bf x})|_{\mu}=\prod_{\mu\in {\cal S}}
|x_1|_{\mu} \le X(n)\ll b_1^{n/2} Y(n),
\end{equation}
while
\begin{equation}
\label{eq:prodge3}
\prod_{\mu\in {\cal S}} |L_{i,\mu}({\bf x})|_{\mu}=\prod_{\mu\in {\cal S}}
|x_i|_{\mu}\le Y(n)\quad \quad {\rm for}~i=3,\dots,M.
\end{equation}
Finally, when $i=2$, and $D=1$, we have, by inequality ~\eqref{ineq:DP2}, that 
\begin{equation}
\label{eq:prod2}
\prod_{\substack{\mu\in {\cal S}\\ \mu< \infty}} |L_{2,\mu}({\bf x})|_{\mu}
\cdot |L_{2,\infty}({\bf x})|_{\infty}\le \frac{1}{(b_1')^{n}}\cdot 
\frac {b_1^{(k-1)n}}{Y(n)}\le \frac{1}{b_1^{n}},
\end{equation}
because $b_1'=b_1^k$, 
while when $i=2$ and $D=2$, we have, again by inequality ~\eqref{ineq:DP2}, 
that 
$$
\prod_{\substack{\mu\in {\cal S}\\ \mu<\infty}} |L_{2,\mu}({\bf x})|_{\mu}
\cdot |L_{2,\infty_1}({\bf x})|_{\infty_1} \cdot |L_{2,\infty_2}({\bf x})|
_{\infty_2}
$$
$$
\le\frac{1}{b_1^{kn}}\cdot 
\(\frac {b_1^{(k-1)n}}{Y(n)}\)^{1/2}\cdot \(\frac {b_1^{(k-1)n}}{Y(n)}\)^{1/2}
\le \frac{1}{b_1^{n}},
$$
which is again inequality ~\eqref{eq:prod2} but for the case $D=2$. 
Inequality ~\eqref{eq:pell} follows now easily by multiplying inequalities 
~\eqref{eq:prod1}, ~\eqref{eq:prodge3} and ~\eqref{eq:prod2}.
We now choose $C_2<\log b_1/(4(M-1))$, and conclude that if $Y(n)<\exp(C_2n)$, 
then $Y(n)^{M-1}\le b_1^{n/4},$ and therefore inequality ~\eqref{eq:pell} 
implies that the inequality
\begin{equation}
\label{eq:pell1}
\prod_{\mu\in {\cal S}} \prod_{i=1}^M |L_{i,\mu}({\bf x})|_{\mu}\ll
\frac{1}{b_1^{n/4}}
\end{equation}
holds. Assuming that $C_2<k\log b_1$, we get that 
$$
{\max}\{|x_i|~|~i=1,\dots,M\}
\ll b_1^{kn}Y(n)\le b_1^{2kn}=(b_1^{n/4})^{8k}.
$$
It follows easily that the above inequality ~\eqref{eq:pell1} implies
that Lemma ~\ref{lem:S-unit} holds for our field $\L$, points ${\bf x}$,  
set of valuations ${\cal S}$ and forms $L_{i,\mu}$ for $i=1,\dots,M$, 
and $\mu\in {\cal S}$, with $\varepsilon=1/(8k+1)$ for all but finitely 
many positive integers $n$.
The conclusion of Lemma ~\ref{lem:S-unit} is that there exist 
finitely many proper subspaces of $\Q^M$ of equations $\sum_{i=1}^M c_ix_i=0$,
with not all the coefficients $c_i$ being zero, and such that all points 
${\bf x}$ satisfying the above inequality ~\eqref{eq:pell1} 
belong to one of these subspaces. 

Assume now that ${\bf x}$ is one of these subspaces of equation 
$\sum_{i=1}^Mc_ix_i=0$. Suppose 
first that $c_1=0$. We then get the equation
$$
\sum_{i=2}^M c_i(b_{i-1}')^n=0,
$$
which gives the set of zeros of
a linear recurrence sequence having a dominant root (note that at least 
one $c_i$ for $i\ge 2$ is nonzero), and as such it can  
have only finitely many positive integer solutions $n$. 

Assume now that $c_1\ne 0$. In this case, we get that 
$$
X(n)=\frac{f_2(n)}{b_1^{(k-1/2)n}}Y(n),
$$
where $f_2(n)\in \Q{\cal E}_{\Z}$ is the form given by 
$$
f_2(n)=-\sum_{i=2}^M c_ic_1^{-1} (b_{i-1}')^n.
$$
Thus, if we write $b=b_1^{(k-1/2)}$, then 
$f_2(n)/b^n=X(n)/Y(n)$. Assume that $C_2<C_1$, where $C_1$ is the 
constant appearing in Lemma ~\ref{lem:CZ1}. Then, 
if $Y(n)<\exp(C_2n)$, and if the above equation has infinitely many 
positive integer solutions $n$, 
it follows, by Lemma ~\ref{lem:fundunit}, 
that $b$ divides every root of $f_2(n)$. In particular, 
we get that $Y(n)$ is bounded. Since we are assuming that this is so for 
infinitely many values of $n$, it follows that there exists a constant value 
$A$ such that $Y(n)=A$ holds for infinitely many values of the positive 
integer $n$. Since the inequality $|X(n)^2-f(n)Y(n)^2|<C$ also holds for all 
these positive integers $n$, it follows that there exists a fixed integer 
$B$ such that both relations 
$X(n)^2-f(n)Y(n)^2=B$ and $Y(n)=A$ hold. In particular,
we conclude that the diophantine equation $f(n)=x^2-B/A^2$ admits infinitely
many solutions $(n,x)$, 
with a positive integer $n$, and a rational number $x$ (namely, 
all the pairs $(n,x)=(n,X(n)/A)$). 
Theorem 3 from ~\cite{CZ1} tells us, in particular, that $f(n)$ 
does not satisfy Hypothesis (H), which is a contradiction.  

The above argument does show that if we choose $C_2$ to be sufficiently small, 
then indeed, for every fixed value of the positive real number $C$, 
all positive integer solutions $(X(n),Y(n))$ 
of the inequality ~\eqref{eq:C} have $Y(n)>\exp(C_2n)$
for all but finitely many values of $n$. 
\end{proof}

\medskip 

\noindent {\bf Remark.} 
It is easy to see that Lemma ~\ref{lem:fundunit} remains 
true even for forms $f(n)\in {\Q}{\cal E}_{\Z}$ satisfying a weaker 
hypothesis then Hypothesis (H), namely that there do not exist 
$j\in \{0,1\}$, $h(n)\in {\Q}{\cal E}_{\Z}$ and $\lambda\in {\Q}$ such that 
$f(2n+j)=h(n)^2+\lambda$ holds identically for all positive integers $n$.

\medskip 

Assume now that $f(n)\in {\Q}{\cal E}_{\Z}$ satisfies Hypothesis (H). 
For every positive integer 
$n$, we write ${\sqrt {f(n)}}=[a_0(n),\dots,a_j(n),\dots]$ 
for the continued fraction expansion 
of ${\sqrt {f(n)}}$. We also write $p_j(n)/q_j(n)$ 
for the $j$th convergent of ${\sqrt {f(n)}}$. The next Lemma is the key
ingredient of the proof of our Theorem ~\ref{thm:Main}, 
as it will show that the first ``sufficiently many'' partial quotients 
$a_j(n)$ are ``small'' for all but finitely many positive integers $n$.

\begin{lem}
\label{lem:partquot}
Let $f(n)\in {\Q}{\cal E}_{\Z}$ be a form satisfying Hypothesis (H). Then 
there exist positive computable 
positive constants $C_3<1$ and $C_4\ge 2$ depending only on $f(n)$, 
such that the following holds.

Assume that $\varepsilon\in (0,C_3)$ is fixed. 

\begin{itemize}
\item[(i)] If $q_j(n)<\exp(C_3\varepsilon n)$, then the inequality 
\begin{equation}
\label{eq:qksmall}
\left|{\sqrt {f(n)}}-\frac{p_j(n)}{q_j(n)}\right|\ge \frac{1}{q_j^2(n)
\exp(\varepsilon n)}
\end{equation}
holds with at most finitely many exceptions in the positive integer 
$n$ (depending on $\varepsilon$).
\item[(ii)] If $\exp(C_3\varepsilon n)\le q_j(n)<\exp(C_3n)$, 
then the inequality 
\begin{equation}
\label{eq:qklarge}
\left|{\sqrt {f(n)}}-\frac{p_j(n)}{q_j(n)}\right|\ge \frac{1}{q_j(n)^{C_4}}
\end{equation}
holds with at most finitely many exceptions in the positive integer $n$
(depending on $\varepsilon$).
\end{itemize}
\end{lem}

\begin{proof} We will deal with both inequalities ~\eqref{eq:qksmall} and 
~\eqref{eq:qklarge} simultaneously. We write 
$Q_j(n)=q_j(n)\exp(\varepsilon n)$ 
in case (i) and $Q_j(n)=q_j(n)^{C_4-1}$ in case (ii). With the notations from
Lemma ~\ref{lem:fundunit}, we have the approximation
$$
{\sqrt {f(n)}}=\alpha\frac{f_1(n)}{b_1^{(k-1/2)n}}+O(b_1^{-n/2}\beta^{-n}).
$$
Thus, the inequality 
$$
\left|{\sqrt {f(n)}}-\frac{p_j(n)}{q_j(n)}\right|<\frac{1}{q_j(n)Q_j(n)}
$$
leads to the inequality
\begin{equation}
\label{eq:approx}
\left|\alpha\frac{f_1(n)}{b_1^{(k-1/2)n}}-
\frac{p_j(n)}{q_j(n)}\right|\ll\frac{1}{q_j(n)Q_j(n)},
\end{equation}
provided that the inequality $(b_1^{1/2}\beta)^n>q_j(n)Q_j(n)$ holds. 
In case (i) this last 
inequality is satisfied if $(C_3+1)\varepsilon<\log(b_1^{1/2}\beta)$, while 
in case (ii) this last 
inequality is satisfied if $C_3C_4<\log(b_1^{1/2}\beta)$. 
Since $\varepsilon<C_3<1$, it follows that in the first case the inequality 
is fulfilled if $2C_3<\log(b_1^{1/2}\beta)$, and since $C_4\ge 2$, we see 
that it suffices that the inequality $C_3C_4<\log(b_1^{1/2}\beta)$ holds. 
From ~\eqref{eq:approx}, we get the inequality
\begin{equation}
\label{ineq:DP3}
\left|b_1^{(k-1/2)n}p_j(n)-\alpha f_1(n)q_j(n)\right|
\ll \frac{b_1^{(k-1/2)n}}{Q_j(n)}.
\end{equation}
Comparing ~\eqref{ineq:DP3} with ~\eqref{ineq:DP2}, we see that 
~\eqref{ineq:DP3} is obtained from ~\eqref{ineq:DP2} by replacing 
$X(n)$ and $Y(n)$ by $p_j(n)$ and $q_j(n)$, respectively, and 
the upper bound $b_1^{(k-1)n}/Y(n)$ on ~\eqref{ineq:DP2} by 
the upper bound $b_1^{(k-1/2)n}/Q_j(n)$. 
We now apply again Lemma ~\ref{lem:S-unit} with the same choices 
of field $\L$, set of places ${\cal S}$, 
forms $L_{i,\mu}$, and integer indeterminates vector ${\bf x}$, 
as in the proof of the Lemma ~\ref{lem:fundunit}. 
Inequality ~\eqref{eq:pell} now becomes
\begin{equation}
\label{eq:pell2}
\prod_{\mu\in {\cal S}} \prod_{i=1}^M |L_{i,\mu}({\bf x})|_{\mu}\ll
 \frac{q_j(n)^{M}}{Q_j(n)}.
\end{equation}

In case (i), the right hand side of ~\eqref{eq:pell2} is 
$q_j(n)^{M-1}/\exp(\varepsilon n)$. Imposing that $C_3<1/(2(M-1))$, 
then $q_j(n)^{M-1}<\exp(\varepsilon n/2)$, and therefore the above 
inequality ~\eqref{eq:pell2} becomes 
\begin{equation}
\label{eq:pellqksmall1}
\prod_{\mu\in {\cal S}} \prod_{i=1}^M |L_{i,\mu}({\bf x})|_{\mu}\ll
 \frac{1}{\exp(\varepsilon n/2)}.
\end{equation}
Assume that $\varepsilon$ is such that 
$C_3\varepsilon <k\log b_1$. Since $\varepsilon<C_3$, it suffices that 
$C_3^2<k\log b_1$. In this case, since $q_j(n)<\exp(C_3\varepsilon n)<
b_1^{kn}$, we get that 
\begin{equation} 
\label{eq:pellqksmall2}
{\rm max}\{|x_i|~|~i=1,\dots,M\}\ll q_j(n)b_1^{kn}\ll b_1^{2kn}=
\exp(\varepsilon n/2)^{\varepsilon^{-1} C_5},
\end{equation}
where $C_5=4k\log b_1$. Hence, from inequalities 
~\eqref{eq:pellqksmall1} and ~\eqref{eq:pellqksmall2}, we get 
\begin{equation}
\label{eq:pellqksmall}
\prod_{\mu\in {\cal S}} \prod_{i=1}^M |L_{i,\mu}({\bf x})|_{\mu}\ll
 \frac{1}{\exp(\varepsilon n/2)}\ll 
{\rm max}\{|x_i|~|~i=1,\dots,M\}^{-\varepsilon C_6},
\end{equation}
where $C_6=C_5^{-1}$.

In case (ii), we may choose 
$C_4=M+2$, and then inequality ~\eqref{eq:pell2} becomes
\begin{equation}
\label{eq:qklarge12}
\prod_{\mu\in {\cal S}} \prod_{i=1}^M 
|L_{i,\mu}({\bf x})|_{\mu}\ll \frac{1}{q_j(n)}\le \frac{1}
{\exp(C_3\varepsilon n)}.
\end{equation}
Assuming that $C_3<k\log b_1$, and that $q_j(n)<\exp(C_3n)$, we note that 
$$
{\rm max}\{|x_i|~|~i=1,\dots,
M\}\ll b_1^{kn}q_j(n)\le b_1^{2kn}=
\exp(C_3\varepsilon n)^{\varepsilon^{-1} C_7},
$$
where $C_7=(2k\log b_1)/C_3$. Thus, inequality ~\eqref{eq:qklarge12} 
implies that the inequality 
\begin{equation}
\label{eq:qklarge1}
\prod_{\mu\in {\cal S}} \prod_{i=1}^M 
|L_{i,\mu}({\bf x})|_{\mu}\ll \frac{1}
{\exp(C_3\varepsilon n)}\ll {\rm max}\{|x_i|~|~i=1,\dots,M
\}^{-\varepsilon C_8} 
\end{equation}
holds with $C_8=C_7^{-1}$.

In either one of the two cases (i) or (ii) we may apply Lemma 
~\ref{lem:S-unit}, and derive that there exist only finitely many 
subspaces of $\Q^M$ of equations $\sum_{i=1}^M c_ix_i=0$, and not
all the coefficients 
$c_i$ being zero, and such that every point ${\bf x}\in \Z^M$ 
satisfying either inequality ~\eqref{eq:pellqksmall} or ~\eqref{eq:qklarge1} 
lies on one of these subspaces. 
 Consider 
now the subspace of equation $\sum_{i=1}^M c_ix_i=0$.
If $c_1=0$, we then get the equation $\sum_{i=2}^{M} c_i(b_{i-1}')^n=0$, 
which has only finitely many positive integer solutions 
$n$ because at least one of the coefficients $c_i$ is nonzero for 
$i=2,\dots,M$. Assume now that $c_1\ne 0$. In this case, we get that 
$$
\frac{p_j(n)}{q_j(n)}=\frac{f_2(n)}{b^n},
$$
where $b=b_1^{(k-1/2)}$, and $f_2(n)=\sum_{i=2}^{M} c_ic_1^{-1} 
(b_{i-1}')^n.$ Assuming that $C_3<C_1$, where $C_1$ appears in Lemma 
~\ref{lem:CZ1}, 
it follows that either the above equation can have only finitely 
many positive integer solutions $n$, or the above equation 
has infinitely many positive integer solutions $n$. In this last case, 
$q_j(n)$ is bounded for all such $n$ and thus,
for large $n$, we are in case (i). It now follows 
that there exists a constant $A$ such that 
$q_j(n)=A$ holds for infinitely many $n$, and we are therefore led 
to the conclusion that the inequality
$$
\left|{\sqrt {f(n)}}-\frac{p_j(n)}{A}\right|\ll \frac{1}{\exp(\varepsilon n)}
$$
holds for infinitely many positive integers $n$. Theorem 3 from ~\cite{CZ1} 
tells us that $f(n)$ does not satisfy Hypothesis (H), which is 
the final contradiction.
\end{proof}

We can now prove our Theorem ~\ref{thm:Main}.

\medskip 

\noindent {\it Proof of Theorem ~\ref{thm:Main}}. We assume again that $b_1$ 
is a perfect square. We write $C_1,~C_2,~C_3,~C_4$ for the constants 
appearing in the statements of Lemmas ~\ref{lem:CZ1}, ~\ref{lem:fundunit} 
and ~\ref{lem:partquot}, respectively. We first note that if 
${\sqrt {f(n)}}$ is a rational number for infinitely many values 
of the positive integer $n$, it follows, by Theorem 3 from ~\cite{CZ1}, 
that there exists a form $h(n)\in {\Q}{\cal E}_{\Z}$ such that 
$f(n)=h(n)^2$. In particular, $f(n)$ does not satisfy Hypothesis (H). 
Assume now that $f(n)$ is a not a square of an integer. In this case, 
${\sqrt {f(n)}}=[a_0(n),\{ a_1(n),\dots,a_{r(n)-1},2a_0(n)\}]$. 
Assume that $r(n)$ does not tend to infinity. Then
there exists a fixed positive integer $r$ such that $r=r(n)$ holds 
for infinitely many positive integers $n$. It is known 
that $p_{r-1}(n)/q_{r-1}(n)$ gives the fundamental unit in the quadratic 
order ${\Q}[{\sqrt {f(n)}}]$. In particular, we have the equation 
$p_{r-1}(n)^2-f(n)q_{r-1}(n)^2=\pm 1$. By Lemma ~\ref{lem:fundunit} with 
$C=1$, it follows that infinitely many positive integers 
$n$ exist such that $q_{r-1}(n)>\exp(C_2n)$. Let $\varepsilon$ 
be a very small number in the interval $(0,C_3)$ to be chosen later. 
By Lemma ~\ref{lem:partquot}, both 
inequalities ~\eqref{eq:qksmall} 
and ~\eqref{eq:qklarge} hold for infinitely many positive integers 
$n$, and for all non negative integers $j$. 
Let $m\le r-1$ be the largest index 
such that the inequality $q_m(n)<\exp(C_3\varepsilon n)$ holds. 
In this case, by inequality 
~\eqref{eq:qksmall}, we get that $q_{m+1}(n)\le q_m(n)\exp(\varepsilon n)<
\exp((C_3+1)\varepsilon n)$, but by the definition of $m$, we also have 
$q_{m+1}(n)\ge \exp(C_3\varepsilon n)$. 
By inequality ~\eqref{eq:qklarge}, we get that $q_{m+2}\le q_{m+1}^{C_4}$ 
once $\varepsilon$ 
is sufficiently small, and, in general, that the inequality
$q_{m+s+1}(n)\le q_{m+s}(n)^{C_4}$ holds provided that 
$q_{m+s}(n)<\exp(C_3 n)$. Assuming therefore that $q_{m+s}(n)<\exp(C_3 n)$, 
we get that $q_{m+s+1}(n)\le q_{m+1}(n)^{C_4^{s}}\le 
\exp(C_4^{s}(C_3+1)\varepsilon n).$ 
Taking $s=r-1$, we get that the inequality
$$
q_r(n)\le q_{m+(r-1)+1}(n)\le \exp(C_4^{r-1}(C_3+1)\varepsilon n)
$$
holds, provided that $C_4^{r-1}(C_3+1)\varepsilon<C_3$. Thus, it suffices 
to choose $\varepsilon$ such that this last inequality is fulfilled. 
However, we also know that $q_r(n)>q_{r-1}(n)\ge \exp(C_2 n)$. 
Hence, if we choose 
$\varepsilon$ such that the inequality $C_4^{r-1}(C_3+1)\varepsilon<C_2$
holds as well, we then obtain a contradiction. 

Thus, $r(n)$ tends to infinity and Theorem ~\ref{thm:Main} is therefore proved.
\qed

\section{Comments and Remarks}

We do not know whether Hypothesis (H) is needed, although 
it is clear that some assumption is necessary in order 
to get the conclusion of Theorem ~\ref{thm:Main}.

Indeed, it is easily checked that for $v(n),~w(n)\in {\cal E}_{\Z}$ and 
$f(n)=v(n)^2w(n)^2+2w(n)$, we have that the relation
$$
{\sqrt {f(n)}}=[v(n)w(n), \{v(n), 2v(n)w(n)\}]
$$
holds for all sufficiently large positive integers $n$.

We are also unable to decide whether for any form $f(n)\in {\cal E}_{\Z}$ 
such that the length $r(n)$ of the period of the continued fraction 
expansion of ${\sqrt {f(n)}}$ remains bounded for infinitely many $n$, 
there must exist $j\in \{0,1\}$ and $f_0(n),\dots,f_{r-1}(n)\in 
{\Q}{\cal E}_{\Z}$ such that the relation
$$
{\sqrt {f(2n+j)}}=[f_0(n),\{f_1(n),\dots,f_{r-1}(n),2f_0(n)\}]
$$
holds for all sufficiently large positive integers $n$.

In the literature, there exist explicit versions of 
Lemma ~\ref{lem:S-unit} (see, for example, ~\cite{Ev},~\cite{EverSchl}), 
which bound the number of possible subspaces occurring 
in Lemma ~\ref{lem:S-unit}. 
Usually, such a bound is of the form $C_5\delta^{-C_6}$. The constant 
$C_6$ depends only on the number of indeterminates $M$, and the number
of places $\#S$, while the constant $C_5$ depends also on the {\it heights} 
of the linear forms $L_{i,\mu}$ for $i=1,\dots,M$, and $\mu\in {\cal S}$. 

It is likely that one could use such results instead of the present 
formulation  of
Lemma ~\ref{lem:S-unit}, in conjunction with upper bounds for the 
zero multiplicities of linearly recurrent sequences (such as 
the results from ~\cite{vdPSchl} and ~\cite{SchlSchm}), to get 
that there exists a function $g(X)$ tending to infinity with $X$, such 
that if $f(n)\in {\cal E}_{\Z}$ satisfies Hypothesis (H) and if $X$ 
tends to infinity, then 
$r(n)\gg g(X)$ holds for all positive integers $n<X$ with $o(X)$ exceptions. 
We point out that a result establishing a lower bound for the 
exponent of the group ${\bf E}(\F_{q^n})$ of points on an elliptic curve 
${\bf E}$ defined over the finite field with $q$ elements $\F_q$, 
and valid for almost all $n$, has been recently established in 
~\cite{LuSh} by a method similar to the one described above.  

Unfortunately, we could not obtain such a result in the present context.

\end{document}